\documentclass[12pt, english]{article}
\usepackage[T2A]{fontenc}
\usepackage[cp1251]{inputenc}
\usepackage{babel}
\usepackage[centertags]{amsmath}
\usepackage{amsfonts}
\usepackage{amssymb}
\usepackage{amsthm}
\usepackage{newlfont}
\usepackage[active]{srcltx}

\newtheorem{theorem}{\bf Theorem}
\newtheorem{lemma}{\bf Lemma}

\begin{document}

\title{On large deviations for combinatorial sums }

\author{Andrei N. Frolov
\footnote{This investigation was supported by RFBR, research project No. 18--01--00393}
\\ Dept. of Mathematics and Mechanics
\\ St.~Petersburg State University
\\ St. Petersburg, Russia
\\ E-mail address: Andrei.Frolov@pobox.spbu.ru}

\maketitle

\abstract{ We investigate asymptotic behaviour of probabilities of large deviations
for normalized combinatorial sums. We find a zone in which  these probabilities
are equivalent to the tail of the standard normal law. Our conditions are
similar to the classical Bernstein condition. The range of the zone of the normal
convergence can be of power order.  }



\section{Introduction}

Let $\{ \left\| X_{nij} \right\|, 1 \leqslant i, j \leqslant n, n=2, 3, \ldots\}$
be a sequence of matrices of independent random variables and
$\{\vec{\pi}_n = (\pi_n(1), \pi_n(2), \ldots, \pi_n(n))$, $n=2, 3, \ldots\}$
be a sequence of random permutations of numbers $1, 2, \ldots,n$.
Assume that $\vec{\pi}_n$ has the uniform distribution on the set of
permutations of $1, 2, \ldots,n$ and it is independent with $\left\| X_{nij} \right\|$ for all $n$. 
Define the combinatorial sum $S_n$ by relation
$$ S_n = \sum\limits_{i=1}^n X_{n i \pi_n(i)}
$$
Under certain conditions, a sequence of distributions of combinatorial sums
converges weakly to the standard normal law. Every such result is called a
combinatorial central limit theorem (CLT).

Investigations in this direction have a long history. 
One can find results on combinatorial CLT in Wald and Wolfowitz [1], 
Noether [2], Hoeffding [3], Motoo [4], Kolchin and Chistyakov [5]. 
Further, non-asymptotic Esseen type bounds have been derived for  
accuracy of normal approximation of distributions of combinatorial sums.
Such results have been obtained in Bolthausen [6], von Bahr [7], 
Ho and Chen [8], Goldstein [9], Neammanee and Suntornchost [10], 
Neammanee and Rattanawong [11], Chen, Goldstein and Shao [12], 
Chen and Fang [13], Frolov [14, 15], and in Frolov [16]
for random combinatorial sums.

Note that if $X_{nij}$ are identically distributed for all
$1\leqslant j \leqslant n$ and $n$, then the combinatorial sum
has the same distribution as that of independent random variables.
This case is well investigated, but one has to take it into account
for  estimation of optimality of derived results.

Besides some partial cases, combinatorial sums have not independent increments.
Hence, it is difficult to use classical methods of proofs for Esseen type inequalities
those are based on bounds for differences of characteristic functions (c.f.).
One usually applies the Stein method. For combinatorial sums, it yields
Esseen type inequalities for random variables with finite third moments.
Applying of the truncation techniques, Frolov [14, 15] 
derived generalizations of these results to the case of finite
moments of order $2+\delta$ and for infinite variations as well.

Every bound in CLT similar to the Esseen inequality yields 
results on asymptotic behaviour for large deviations coinciding
with that for tail of the normal law in a logarithmic zone.
Such results are usually called moderate deviations.
Moderate deviations for combinatorial sums have been investigated in Frolov [17].

In this paper, we derive new results on the asymptotic behaviour for large deviations
of combinatorial sums in power zones. Note that ranges of power zones are 
powers from some characteristic similar to the Lyapunov ratio. Indeed, we deal with
non-identically distributed random variables. Even for sums of independent
random variables, ranges of zones of the normal convergence  depend on the Lyapunov ratios.
For identically distributed random variables, this yields that the ranges are powers 
from the number of summands. But the last case corresponds to the classical theory
for sums of independent random variables and it is not new therefore.

In our proofs, we will use the method of conjugate distributions.
Note that von Bahr [7] developed a method to bound distances between
c.f.'s of normalized combinatorial sums and normal law.
Assuming that random variables are bounded or satisfy certain 
analogue of the classical Bernstein condition, we conclude that
moment generating functions (m.g.f.) of normalised combinatorial sums
are analytic in a circle of the complex plane. Adopting the Bahr's method,
we will bound the difference between m.g.f.'s in some circle. 
In view of the analytic property, this will also give bounds for derivatives
of m.g.f.'s. Hence, we will arrive at desired asymptotics for
m.g.f.'s and their first and second logarithmic derivatives
which are means and variations of random variables being conjugate for normalized combinatorial sums.
Then we will estimate a closeness of distributions of conjugate random variables and 
the standard normal law. Using relationship between distributions and conjugate ones, 
we will derive the asymptotics of large deviations under consideration.

\section{Results}

Let $\{ \left\|X_{nij}\right\|, 1 \leqslant i, j \leqslant n, n=2,3,\ldots\}$ be 
a sequence of matrices of independent random variables such that
\begin{eqnarray}\label{10}
\sum\limits_{i=1}^n \mathbf{E} X_{nij} = \sum\limits_{j=1}^n \mathbf{E} X_{nij}=0
\end{eqnarray}
for all $n$. Let 
$\{\vec{\pi}_n = (\pi_n(1), \pi_n(2), \ldots, \pi_n(n))$, $n=2, 3, \ldots\}$
be a sequence of random permutations of numbers $1,2, \ldots,n$.
Assume that $\vec{\pi}_n$ has the uniform distribution on the set of permutation
$P_n$ and it is independent with $\left\|X_{nij}\right\|$ for all $n$. 
Put
$$ S_n = \sum\limits_{i=1}^n X_{n i \pi_n(i)}.
$$
It is not difficult to check that
$$ \mathbf{E} S_n=0, \quad \mathbf{D} S_n = \mathbf{E} S_n^2 - (\mathbf{E} S_n)^2 =
\frac{1}{n-1}\sum\limits_{i,j=1}^n (\mathbf{E} X_{nij})^2 + 
\frac{1}{n} \sum\limits_{i,j=1}^n \mathbf{D} X_{nij} . $$
Hence, condition (\ref{10}) yields that combinatorial sums are centered at zero.
Moreover,
$$ \mathbf{D} S_n= 
\frac{1}{n(n-1)}\sum\limits_{i,j=1}^n (\mathbf{E} X_{nij})^2 + \frac{1}{n} \sum\limits_{i,j=1}^n \mathbf{E} X_{nij}^2. $$
If $\mathbf{D} S_n \to\infty$ as $n\to \infty$, then the main part of the variance is the
normalized sum of second moments
$$
B_n = \frac{1}{n} \sum\limits_{i,j=1}^n \mathbf{E} X_{nij}^2.
$$
Therefore, in the sequel, we will use $\{B_n\}$ as norming sequence for
$S_n$.

Our main result is as follows. 

\begin{theorem}\label{t1}
Let $ \{M_n\} $ be a non-decreasing sequence of positive numbers  such that
for $s=1,2,3$, inequalities
\begin{eqnarray}\label{bern}
\left| \mathbf{E} X_{nij}^k \right| \leqslant D k! M_n^{k-s} \mathbf{E} |X_{nij}|^s
\end{eqnarray}
hold for all  $k \geqslant s$, $1 \leqslant i, j \leqslant n$ and  $n \geqslant 2$,
where $D$ is an absolute positive constant.
Put
$$ \gamma_n =
\max\left\{ \max_{i,j} \frac{\sqrt{n}}{\sqrt{B_n}} \mathbf{E}|X_{nij}|,\;
\max_i  \sum\limits_{j=1}^n  \frac{\mathbf{E} X_{nij}^2}{B_n},\;
\max_j \sum\limits_{j=1}^n  \frac{\mathbf{E} X_{nij}^2}{B_n},\;
\sum\limits_{i,j=1}^n  \frac{\mathbf{E} |X_{nij}|^3}{\sqrt{n} B_n^{3/2}} 
\right\}.
$$

Then for every sequence of positive numbers $\{u_n\}$ with
$u_n \to \infty$, $u_n^3=o(\sqrt{n}/\gamma_n)$ and $u_n=o(\sqrt{B_n}/M_n)$
as $n\to\infty$, relation
\begin{eqnarray}\label{20}
\mathbf{P}\left(S_n \geqslant u_n\sqrt{B_n}\right)  \sim 1-\Phi(u_n)
\quad\mbox{as}\quad n \to \infty,
\end{eqnarray}
holds, where $\Phi(x)$ is the standard normal distribution function.
\end{theorem}

Note that $\gamma_n \geqslant 1$. This follows from the inequality
$\max\limits_i  \sum\limits_{j=1}^n  \mathbf{E} X_{nij}^2 \geqslant B_n$.
Indeed, assuming that
$\max\limits_i  \sum\limits_{j=1}^n  \mathbf{E} X_{nij}^2 < B_n$, we arrive at
the incorrect inequality
$\sum\limits_{i,j=1}^n  \mathbf{E} X_{nij}^2 < n B_n = \sum\limits_{i,j=1}^n  \mathbf{E} X_{nij}^2$.

Bahr [7] proved the following Esseen type inequality:
$$ \sup_{x} \left| \mathbf{P}\left(S_n < x \sqrt{B_n}\right) - \Phi(x)
\right| \leqslant A \frac{\gamma_n}{\sqrt{n}},
$$
where $A$ is an absolute positive constant. Hence the condition
$u_n^3=o(\sqrt{n}/\gamma_n)$ as $n\to\infty$ is natural for relation (\ref{20}),
giving exact (non-logarithmic) asymptotics of large deviations.
For identically distributed $X_{nij}$, this condition turns to 
$u_n = o(n^{1/6})$ as $n\to\infty$.

Note that the conditions $u_n^3=o(\sqrt{n}/\gamma_n)$ and $u_n \to \infty$ as $n\to\infty$
imply $\gamma_n/\sqrt{n} \to 0$ as $n\to\infty$.

Theorem \ref{t1} is stronger than the results in Frolov [14] 
since the zone of normal convergence may be of power order while
it is logarithmic in [14]. Of course, this requires stronger moment assumptions.

Condition (\ref{bern}) is an analogue of the Bernstein condition which
is a form of existence for the exponential moment. In classical theory,
one mainly deals with centered random variables and the Berstein condition
yields that the logarithm of the m.g.f. is asymptotically a quadratic function
at zero. For combinatorial CLT, it is principally important that summands could be non-centered 
and even degenerate sometimes. In this case, the logarithm of m.g.f. 
may be a linear function in a neighbourhood of zero provided the mean is not zero.
 
One can rewrite inequalities (\ref{bern}) for $k \geqslant 3$ as follows: 
\begin{eqnarray*}
\left| \mathbf{E} X_{nij}^k \right| \leqslant D k! M_n^{k}
\min_{1\leqslant s \leqslant 3} \frac{\mathbf{E} |X_{nij}|^s}{M_n^s}.
\end{eqnarray*}
Hence, the Lyapunov inequality implies that the next condition is
sufficient for (\ref{bern}): the inequalities 
$\mathbf{E} X_{nij}^2 \leqslant 2 D M_n \mathbf{E} |X_{nij}|$ 
and
\begin{eqnarray*}
\left| \mathbf{E} X_{nij}^k \right| \leqslant D k! M_n^{k} 
\min\left\{
\frac{\mathbf{E} |X_{nij}|}{M_n},
\left(\frac{\mathbf{E} |X_{nij}|}{M_n}\right)^3
\right\}
\end{eqnarray*}
hold for all  $k \geqslant 3$,  $1 \leqslant i, j \leqslant n$ and $n \geqslant 2$.

Consider two important examples in which condition (\ref{bern}) is satisfied.

{\bf 1. Bounded random variables.} If there exists a non-decreasing sequence 
of positive constants $\{M_n\}$ such that $\mathbf{P}(|X_{nij}| \leqslant M_n)=1$  
for all $1 \leqslant i, j \leqslant n$ and $n \geqslant 2$, 
then condition (\ref{bern}) holds. For degenerate case with $\mathbf{P}(X_{nij}=c_{nij})=1$
for all $1 \leqslant i, j \leqslant n$ and $n \geqslant 2$, condition (\ref{bern})
is fulfilled with $M_n = \max_{i,j} |c_{nij}|$ for every $n$.
 
{\bf 2. Exponential random variables.}
Let $ \xi $ and $ \eta $ be random variables having the exponential distributions
with the parameters $ \alpha $ and $ \beta $ correspondingly.
Assume that each random variable in every matrix $\|X_{nij}\|$ 
has one from four distributions of random variables $ \xi $, $ -\xi $, $ \eta $ and $ -\eta $.
Since $ \mathbf{E} \xi^k = \alpha^{-k}$ and $ \mathbf{E} \eta^k = \beta^{-k}$ for all $k$,
condition (\ref{bern}) holds with с $M_n=1/\min(\alpha,\beta)$.
One can easily expand this example for a larger number of exponential
distributions using for construction of matrices of X's.
Parameters of these distributions may depend on  $n$. 
Moreover, one can easily replace exponential distributions by Gamma ones.

Note that $ \gamma_n $ has an order of
$\max\{\sqrt{n}/B_n, (\sqrt{n}/B_n)^3\}$ in the last example.
It is also clear that the behaviour of $ \gamma_n $ will be similar
when every random variable $X_{nij}$ has one from $k$ given
distributions. In the last case, one says about $k$-sequences
of matrix $\{\|X_{nij}\|\}$.

\section{ Proofs}

For all $i$, $j$ and $n$, put
$$ \varphi_{nij}(z)=\mathbf{E} e^{z X_{nij}},\quad
\varphi_n(z) = \mathbf{E} e^{z \frac{S_n}{\sqrt{B_n}}}, \quad z \in \mathbb{C},
$$
where $\mathbb{C}$ is the set of complex numbers.
We have
\begin{eqnarray}
\label{efi}
e^{-\frac{z^2}{2}} \varphi_n(z) \!=\! \frac{1}{n!} \sum_{p_n \in P_n}
\prod_{i=1}^n \left\{ e^{-\frac{z^2}{2n}} \varphi_{ni p_n(i)}\left(\frac{z}{\sqrt{B_n}}\right)\right\}
\!=\! \frac{1}{n!} \sum_{p_n \in P_n}
\prod_{i=1}^n \left\{ 1+ b_{ni p_n(i)}\right\}.
\end{eqnarray}
Note that the last sum is the permanent of the matrix $\|1+b_{nij}\|$.
To investigate its behaviour we will use the following result.

\begin{lemma}\label{l1}
Let $X$ be a random variable such that for $s=1,2,3$ the inequalities
\begin{eqnarray}\label{bern1}
| \mathbf{E} X^k| \leqslant D k! M^{k-s} \mathbf{E}|X|^s
\end{eqnarray}
hold for all $k \geqslant s$, where $D$ and $M$ are positive constants.

Then $\mathbf{E} e^{u X}$ is an analytic function in the circle $|u| \leqslant 1/(4 M)$ and
for every $u, v \in \mathbb{C}$ with $|v|\leqslant 1/2$ and $|u| \leqslant 1/(8 M)$,
the inequalities
\begin{eqnarray*}
&&
\left| \mathbf{E} e^{u X -\frac{v^2}{2}} -1\right| \leqslant C_1 ( |u| \mathbf{E}|X| + |v|),
\\ &&
\left| \mathbf{E} e^{u X -\frac{v^2}{2}} -1- u \mathbf{E} X \right| 
\leqslant C_2 ( |u|^2 \mathbf{E} X^2 + |v|^2),
\\ &&
\left| \mathbf{E} e^{u X -\frac{v^2}{2}} -1- u \mathbf{E} X + \frac{v^2}{2}
- \frac{u^2}{2} \mathbf{E} X^2
\right| 
\leqslant C_3 ( |u|^3 \mathbf{E} |X|^3 + |v|^3)
\end{eqnarray*}
hold, where constants $C_i$ depends on $D$ and do not depend on $M$.
\end{lemma}

{\bf Proof.} By inequality (\ref{bern1}) and Stirling's formula, we have
\begin{eqnarray*}
E|X|^k \leqslant \sqrt{E X^{2k}} \leqslant \sqrt{D (2k)! M^{2k-1} \mathbf{E}|X|}
\leqslant D_1 k! (2 M)^k
\end{eqnarray*}
for all $k \geqslant 1$, where the constant $D_1$ depends on $D$, $M$ and $\mathbf{E} |X|$.
Hence, the series 
$$ \sum_{k=0}^\infty \frac{|u|^k}{k!} \mathbf{E} |X|^k
$$
converges in the circle $|u|\leqslant 1/(4 M)$.
Put 
$$ e_n(z) = \sum_{k=0}^n \frac{z^k}{k!}.
$$
Then $e_n(|u X|) \uparrow e^{|u X|}$ a.s. The monotone convergence theorem yields that
$$ \mathbf{E}  e^{|u X|} = \lim_{n \to \infty} \mathbf{E} e_n(|u X|) = 
\sum_{k=0}^\infty \frac{|u|^k}{k!} \mathbf{E} |X|^k
$$ 
in the circle $|u|\leqslant 1/(4 M)$. In view of $|e^{u X}| \leqslant  e^{|u X|}$,
the Lebesgue dominate convergence theorem implies that
$$ \mathbf{E}  e^{u X} = \lim_{n \to \infty} \mathbf{E} e_n(u X) = 
\sum_{k=0}^\infty \frac{u^k}{k!} \mathbf{E} X^k
$$
in the circle $|u|\leqslant 1/(4 M)$.

Put $W=u X-v^2/2$. For $s=1,2,3$ and $k\geqslant s$, we have 
\begin{eqnarray*}
\left|\mathbf{E} W^k\right|\! = \!
\left|\mathbf{E} \sum\limits_{j=0}^k C_k^j (uX)^j \left(-\frac{v^2}{2}\right)^{k-j}
\right| 
\!\leqslant\! \sum\limits_{j=0}^k C_k^j |u|^j |\mathbf{E}X^j| \left(\frac{|v|^2}{2}\right)^{k-j} \!=
T'_{ks}+T''_{ks},
\end{eqnarray*}
where
\begin{eqnarray*}
T'_{ks} = \sum\limits_{j=0}^{s-1} C_k^j |u|^j |\mathbf{E}X^j| \left(\frac{|v|^2}{2}\right)^{k-j},
\quad
T''_{ks} = \sum\limits_{j=s}^k C_k^j |u|^j |\mathbf{E}X^j| \left(\frac{|v|^2}{2}\right)^{k-j}.
\end{eqnarray*}
By inequalities (\ref{bern1}), we get
\begin{eqnarray*}
&& 
T''_{ks} \leqslant D k! |u|^s \mathbf{E} |X|^s
\sum\limits_{j=s}^k C_k^j \left|u M\right|^{j-s} \left(\frac{|v|^2}{2}\right)^{k-j}
\\ && 
\leqslant D k! |u|^s \mathbf{E} |X|^s
\sum\limits_{j=0}^{k-s} C_{k}^{j+s} \left|u M\right|^{j} \left(\frac{|v|^2}{2}\right)^{k-j-s}
\\ && 
\leqslant D k! |u|^s \mathbf{E}|X|^s \sum\limits_{j=0}^{k-s} k^s
C_{k-s}^{j} \left|u M\right|^{j} \left(\frac{|v|^2}{2}\right)^{k-j-s}
\\ &&
\leqslant D k!  |u|^s \mathbf{E} |X|^s k^s
\left(\left|u M\right| + \frac{|v|^2}{2} \right)^{k-s}
\leqslant D k!  |u|^s \mathbf{E}|X|^s k^s 4^{-k+s}
\end{eqnarray*}
for all $k \geqslant s$ and $s=1,2,3$.

Since $|v| \leqslant 1/2$, we have
\begin{eqnarray*}
T'_{k1} = \frac{|v|^{2k}}{2^k} \leqslant 2 |v| 8^{-k}
\end{eqnarray*}
for all $k \geqslant 1$.

Hence,
\begin{eqnarray*}
&&
\left| \mathbf{E} e^{u X -\frac{v^2}{2}} -1\right| 
= \left| \sum\limits_{k=1}^\infty \frac{1}{k!} \mathbf{E} W^k\right| \leqslant
\sum\limits_{k=1}^\infty \frac{1}{k!} |\mathbf{E} W^k| 
\leqslant \sum\limits_{k=1}^\infty \frac{1}{k!} (T'_{k1}+T''_{k1})
\\ &&
\leqslant 4 D |u| \mathbf{E}|X| \sum\limits_{k=1}^\infty k 4^{-k} + 2 |v|
\sum\limits_{k=1}^\infty \frac{8^{-k}}{k!}
\leqslant
C_1 (|u| \mathbf{E} |X|+ |v| ).
\end{eqnarray*}
The first inequality follows.

Making use of the inequality $2 a\leqslant a^2+1$ for
$a= |u| |\mathbf{E}X|$, the Lyapunov inequality and $|v|\leqslant 1/2$,
we obtain
\begin{eqnarray*}
&&
T'_{k2} = k |u| |\mathbf{E}X| \frac{|v|^{2k-2}}{2^{k-1}} +\frac{|v|^{2k}}{2^k}
\leqslant
\frac{k}{2} |u|^2 (\mathbf{E}X)^2 \frac{|v|^{2k-2}}{2^{k-1}} +
 \frac{k}{2} \frac{|v|^{2k-2}}{2^{k-1}} +\frac{|v|^{2k}}{2^k}
\\ && 
\leqslant
 4 k 8^{-k} |u|^2 \mathbf{E}X^2  + 
16 k |v|^{2} 8^{-k}
+ 4 |v|^{2} 8^{-k}
\leqslant
 4 k 8^{-k} |u|^2 \mathbf{E}X^2  + 
20 k |v|^{2} 8^{-k}
\end{eqnarray*}
for all $k \geqslant 2$.

It follows that
\begin{eqnarray*}
&&
\left| \mathbf{E} e^{u X -\frac{v^2}{2}} -1- u \mathbf{E} X \right| 
= \left|  -\frac{v^2}{2}+ 
\sum\limits_{k=2}^\infty \frac{1}{k!} \mathbf{E} W^k\right| 
\leqslant \frac{|v|^2}{2}+ \sum\limits_{k=2}^\infty \frac{1}{k!} (T'_{k2}+T''_{k2})
\\ &&
\leqslant \frac{|v|^2}{2}+ 16 (D+1)
|u|^2 \mathbf{E} X^2  \sum_{k=2}^\infty k^2 4^{-k}
+ 20 |v|^2 \sum_{k=2}^\infty \frac{8^{-k}}{k!}
\leqslant C_2
(|u|^2 \mathbf{E}X^2 + |v|^2).
\end{eqnarray*}
The second inequality is proved.

Applying the inequality $2 a\leqslant a^2+1$ for
$a= |u| |\mathbf{E}X|$ and the Lyapunov inequality, we have
\begin{eqnarray*}
&& T'_{k3} = \frac{k(k-1)}{2} |u|^2 \mathbf{E}X^2 \frac{|v|^{2k-4}}{2^{k-2}}
+ k |u| |\mathbf{E}X| \frac{|v|^{2k-2}}{2^{k-1}}
+\frac{|v|^{2k}}{2^k}
\\ &&
\leqslant
\frac{k(k-1)}{2} |u|^2 \mathbf{E}X^2 \frac{|v|^{2k-4}}{2^{k-2}}
+ \frac{k}{2} |u|^2 \mathbf{E}X^2 \frac{|v|^{2k-2}}{2^{k-1}} + \frac{k}{2} 
\frac{|v|^{2k-2}}{2^{k-1}}
+\frac{|v|^{2k}}{2^k}
\\ &&
\leqslant
k^2 |u|^2 \mathbf{E}X^2 \frac{|v|^{2k-2}}{2^{k-1}} + \frac{k}{2} 
\frac{|v|^{2k-2}}{2^{k-1}}
+\frac{|v|^{2k}}{2^k}.
\end{eqnarray*}
Using  $2 a\leqslant a^2+1$ for $a= |u|^2 \mathbf{E}X^2$, 
the Lyapunov inequality, inequality (\ref{bern1}) and $|v|\leqslant 1/2$, we further get
\begin{eqnarray*}
&& T'_{k3}
\leqslant
 \frac{k^2}{2} |u|^4 (\mathbf{E}X^2)^2 \frac{|v|^{2k-2}}{2^{k-1}} 
+ \frac{k^2}{2} \frac{|v|^{2k-2}}{2^{k-1}} +
\frac{k}{2} \frac{|v|^{2k-2}}{2^{k-1}}
+\frac{|v|^{2k}}{2^k}
\\ &&
\leqslant
 \frac{k^2}{2} |u|^4 \mathbf{E}X^4 \frac{|v|^{2k-2}}{2^{k-1}} 
+ 2 k^2 \frac{|v|^{2k-2}}{2^{k-1}} 
\leqslant
\frac{k^2}{2} |u|^4 D 4! M \mathbf{E}|X|^3 \frac{|v|^{2k-2}}{2^{k-1}} 
+ 2^7 k^2 |v|^{3} 8^{-k}
\\ && 
\leqslant  
3 D \frac{k^2}{2} |u|^3 \mathbf{E}|X|^3 \frac{|v|^{2k-2}}{2^{k-1}} 
+ 2^7 k^2 |v|^{3} 8^{-k}
\leqslant
12 D k^2 |u|^3 \mathbf{E}|X|^3 8^{-k} 
+ 2^7 k^2 |v|^{3} 8^{-k}
\end{eqnarray*}
for all $k \geqslant 3$.

It yields that
\begin{eqnarray*}
&&
\left| \mathbf{E} e^{u X -\frac{v^2}{2}} -1- u \mathbf{E} X + \frac{v^2}{2}
- \frac{u^2}{2} \mathbf{E} X^2
\right| 
= \left|  -\frac{u^2}{2} \mathbf{E} X^2 + \frac{1}{2} \mathbf{E} W^2+
\sum\limits_{k=3}^\infty \frac{1}{k!} \mathbf{E} W^k\right| 
\\ &&
\leqslant
\left|  -\frac{u^2}{2} \mathbf{E} X^2 + \frac{1}{2} \mathbf{E} W^2\right| +
\sum_{k=3}^\infty \frac{1}{k!} (T'_{k3}+T''_{k3})
\\ &&
\leqslant
\left|  -\frac{u^2}{2} \mathbf{E} X^2 + \frac{1}{2} \mathbf{E} W^2\right| +
76 D  |u|^3 \mathbf{E}|X|^3 
\sum_{k=3}^\infty k^3 4^{-k} + 2^7 |v|^{3} \sum_{k=3}^\infty
k^2  \frac{8^{-k}}{k!}
\\ &&
\leqslant \frac{|v|^2}{4} |\mathbf{E}(W+uX)|+ 
76 D  |u|^3 \mathbf{E}|X|^3 
\sum_{k=3}^\infty k^3 4^{-k} + 2^7 |v|^{3} \sum_{k=3}^\infty
k^2  \frac{8^{-k}}{k!}
\\ &&
\leqslant
C_3 (|u|^3 \mathbf{E} |X|^3+|v|^3).
\end{eqnarray*}
The lemma is proved. $\Box$

{\bf Proof of Theorem \ref{t1}.}
By Lemma \ref{l1} with $X=X_{nij}$, $u=z/\sqrt{B_n}$, $v=z/\sqrt{n}$ and $M=M_n$,
for all $n$, $i$ and $j$, the inequalities
\begin{eqnarray}
&&
\label{b1}
\left|b_{nij} \right| \leqslant C_1 
\left( \frac{|z|}{\sqrt{B_n}} \mathbf{E}|X_{nij}| + \frac{|z|}{\sqrt{n}}\right),
\\ &&
\label{b2}
\left| b_{nij} - \frac{z}{\sqrt{B_n}} \mathbf{E} X_{nij} \right| 
\leqslant C_2 
\left( \frac{|z|^2}{B_n} \mathbf{E} X_{nij}^2 + \frac{|z|^2}{n}\right),
\\ &&
\label{b3}
\left| b_{nij} - \frac{z}{\sqrt{B_n}} \mathbf{E} X_{nij} + \frac{z^2}{2n}
- \frac{z^2}{2 B_n} \mathbf{E} X_{nij}^2
\right| 
\leqslant C_3 
\left( \frac{|z|^3}{B_n^{3/2}} \mathbf{E}|X_{nij}|^3 + \frac{|z|^3}{n^{3/2}}\right)
\end{eqnarray}
hold for every $z$ in the circle $|z| \leqslant \min\{ \sqrt{n}, \sqrt{B_n}/M_n\}/8$.
Relations (\ref{b2}) and (\ref{10}) imply that
\begin{eqnarray}\label{120}
|b_{n \cdot j}| = \left| \sum\limits_{i=1}^n b_{nij}\right|
= \left| \sum\limits_{i=1}^n \left( b_{nij}- \frac{z}{\sqrt{B_n}} \mathbf{E} X_{nij}\right)
\right| \leqslant C_2 |z|^2 \left(1+
\frac{1}{B_n} \sum\limits_{i=1}^n  \mathbf{E} X_{nij}^2 \right)
\end{eqnarray}
and
\begin{eqnarray}\label{130}
|b_{n i \cdot}| = \left| \sum\limits_{j=1}^n b_{nij}\right|
\leqslant C_2 |z|^2 \left(1+
\frac{1}{B_n} \sum\limits_{j=1}^n  \mathbf{E} X_{nij}^2 \right).
\end{eqnarray}
It follows from relations (\ref{b3}) and (\ref{10}) that
\begin{eqnarray}
&& \nonumber
|b_{n \cdot \cdot}| = \left| \sum\limits_{i,j=1}^n b_{nij}\right|
= \left| \sum\limits_{i,j=1}^n \left( b_{nij}- \frac{z}{\sqrt{B_n}} \mathbf{E} X_{nij}
+ \frac{z^2}{2n}- \frac{z^2}{2 B_n} \mathbf{E} X_{nij}^2 \right)
\right| 
\\ && \label{140}
\leqslant C_3 |z|^3 \sqrt{n} \left(1+
\frac{1}{\sqrt{n} B_n^{3/2}} \sum\limits_{i,j=1}^n  \mathbf{E} |X_{nij}|^3 \right).
\end{eqnarray}

From relations (\ref{b1}) and (\ref{120})---(\ref{140}), the definition of 
$ \gamma_n $ and $ \gamma_n \geqslant 1 $, we have
\begin{eqnarray}
&& \label{150}
|b_{nij}| \leqslant C_1 \left( \gamma_n+1 \right) \frac{|z|}{\sqrt{n}}
\leqslant 2 C_1 \gamma_n \frac{|z|}{\sqrt{n}},
\\ && \label{160}
\sum_{j=1}^n |b_{n\cdot j}| \leqslant C_2 \left( \gamma_n+1 \right) |z|^2 n 
\leqslant 2 C_2 \gamma_n |z|^2 n,
\\ && \label{170}
\sum_{i=1}^n |b_{n i \cdot }| 
\leqslant 2 C_2 \gamma_n |z|^2 n,
\\ && \label{180}
|b_{n \cdot \cdot }| 
\leqslant C_3 \left( \gamma_n+1 \right) |z|^3 \sqrt{n}
\leqslant 2 C_3 \gamma_n |z|^3 \sqrt{n}.
\end{eqnarray}

Note that the function $\varphi_n(it)$, $t\in \mathbb{R}$,  
is the c.f. for the normalized combinatorial sum.
In Bahr [7], relations (\ref{efi}) and (\ref{150})---(\ref{180}) for $z=i t$ and $t \geqslant 0$
have been used to bound the distance between $\varphi_n(i t)$ and the
c.f. of the standard normal law. The bounds for $b_{nij}$ from there will
coincide with our ones provided we change  $t$ by $|z|$.  Hence, we borrow
one further bound from [7] with a formal replacing  $t$ by $|z|$. 
We use the first formula from p. 137 in [7] with  $2 C_3 \gamma_n$ instead of
$\delta$. Then we have
\begin{eqnarray*}\label{190}
\left|e^{-\frac{z^2}{2}} \varphi_n(z)-1 \right|
\leqslant \frac{1}{2} \sum\limits_{k=1}^n
\frac{1}{k!} \left(8 e C_3 \frac{\gamma_n |z|^3}{\sqrt{n}}\right)^k
+ \frac{1}{2} \sum\limits_{k=1}^n 
\left(4 e^2 C_4 \frac{\gamma_n |z|}{\sqrt{n}}\right)^k
\end{eqnarray*}
for all $z$ in the circle $|z| \leqslant \min\{ \sqrt{n}, \sqrt{B_n}/M_n\}/8$,
where $C_4$ is an absolute positive constant. 
Hence, 
\begin{eqnarray*}
\left|e^{-\frac{z^2}{2}} \varphi_n(z)-1 \right|
\leqslant \frac{1}{2} \sum\limits_{k=1}^n
\frac{1}{k!} \left(C_5 \frac{\gamma_n |z|^3}{\sqrt{n}}\right)^k
+ \frac{1}{2} \sum\limits_{k=1}^n 
\left( C_5 \frac{\gamma_n |z|}{\sqrt{n}}\right)^k,
\end{eqnarray*}
where $C_5=\max\{8e C_3, 4 e^2 C_4\}$. If $|z| \leqslant \sqrt{n}/(2 C_5 \gamma_n)$, then
\begin{eqnarray*}
\left|e^{-\frac{z^2}{2}} \varphi_n(z)-1 \right|
\leqslant \frac{1}{2} C_5 \frac{\gamma_n |z|^3}{\sqrt{n}}
\sum\limits_{k=1}^\infty
\frac{1}{(k-1)!} \left(C_5 \frac{\gamma_n |z|^3}{\sqrt{n}}\right)^{k-1}
+  C_5 \frac{\gamma_n |z|}{\sqrt{n}} \sum_{k=1}^\infty 2^{-k}.
\end{eqnarray*}
It follows that
\begin{eqnarray}\label{210}
\left|e^{-\frac{z^2}{2}} \varphi_n(z)-1 \right|
\leqslant 
C_5 \frac{\gamma_n |z|^3}{\sqrt{n}} e^{C_5 \frac{\gamma_n |z|^3}{\sqrt{n}}}
+ C_5 \frac{\gamma_n |z|}{\sqrt{n}}= g_n(|z|)
\end{eqnarray}
for all $z$ in the circle  $|z| \leqslant C_7 \min\{ \sqrt{n}/\gamma_n,
\sqrt{B_n}/M_n\}=y_n$.

Let $\{x_n\}$ be a sequence of positive numbers that will be chosen later.
Assume that $x_n \leqslant y_n/16$. The function $ f_n(z) = e^{-\frac{z^2}{2}} \varphi_n(z)-1 $
is analytic in the circle $|z|\leqslant 16 x_n$. Hence,
$$ f_n(z) = \sum\limits_{k=1}^\infty a_{nk} z^k,\quad
f'_n(z) = \sum\limits_{k=1}^\infty a_{nk} k z^{k-1},\quad
f''_n(z) = \sum\limits_{k=1}^\infty a_{nk} k(k-1) z^{k-2},
$$
where by the Cauchy inequalities, the coefficients $a_{nk}$ satisfy to the relations
$$ |a_{nk}| \leqslant (8 x_n)^{-k} \sup_{|z|=8 x_n} |f(z)| \leqslant (8 x_n)^{-k} g_n(8 x_n).
$$
Put
\begin{eqnarray*}
m_n(z) = \frac{\varphi_n'(z)}{\varphi_n(z)}, \quad
\sigma_n^2(z) = \frac{\varphi_n''(z)}{\varphi_n(z)} - 
\left(\frac{\varphi_n'(z)}{\varphi_n(z)}\right)^2, \quad z \in \mathbb{C}.
\end{eqnarray*}
Then, in the circle $|z|\leqslant 4 x_n$, the inequalities
\begin{eqnarray}
&& \label{220}
\hspace*{-\parindent}
\left| \left(m_n(z)-z\right)
e^{-z^2/2} \varphi_n(z)\right| =
|f'_n(z)| \leqslant \frac{g_n(8x_n)}{4 x_n} 
\sum\limits_{k=1}^\infty k 2^{-k}
\leqslant C_8  \frac{g_n(8 x_n)}{ x_{n}},
\\ && \nonumber
\hspace*{-\parindent}
\left| \left(\sigma_n^2(z)-1+(m_n(z)-z)^2\right)
e^{-z^2/2} \varphi_n(z)\right| =
|f''_n(z)| \leqslant 
\frac{g_n(8 x_n)}{(4 x_n)^2} 
\sum\limits_{k=1}^\infty k(k-1) 2^{-k}
\\ && \label{230}
\hspace*{-\parindent}
\leqslant C_9  \frac{g_n(8 x_n)}{ x_{n}^2}
\end{eqnarray}
hold. This and inequality (\ref{210}) yield that
\begin{eqnarray} 
&& \nonumber
\left|m_n(z)-z\right|
\leqslant 
\left| \left(m_n(z)-z\right)
e^{-z^2/2} \varphi_n(z)\right| +
\left| \left(m_n(z)-z\right)
\left(e^{-z^2/2} \varphi_n(z)-1\right)\right|
\\ &&
\nonumber \leqslant
C_8  \frac{g_n(8 x_n)}{ x_{n}} +
\left|m_n(z)-z\right|
g_n(8 x_n)
\end{eqnarray}
for $|z|\leqslant 4 x_n$.
Hence,
\begin{eqnarray}\label{240}
&& 
\left|m_n(z)-z\right|
\leqslant 
C_8  \frac{g_n(8x_n)}{ x_{n} (1-g_n(8x_n))}
\end{eqnarray}
for $|z|\leqslant 4 x_n$.

Further, making use of relations (\ref{220})--(\ref{240}), we get
\begin{eqnarray*}
&& \nonumber
\left| \sigma_n^2(z)-1\right| \leqslant
\left| \left(\sigma_n^2(z)-1\right)
\left(e^{-z^2/2} \varphi_n(z)-1\right)\right|+
\\ && \nonumber
\left| \left(\sigma_n^2(z)-1+(m_n(z)-z)^2\right)
e^{-z^2/2} \varphi_n(z)\right|
+
\left| \left(m_n(z)-z\right)^2
e^{-z^2/2} \varphi_n(z)\right|
\\ &&
\leqslant
\left|\sigma_n^2(z)-1\right| g_n(2x_n) +
C_9  \frac{g_n(8x_n)}{ x_{n}^2} +
C_8^2   \frac{g_n^2(8x_n)}{ x_{n}^2 (1-g_n(8x_n))}
\end{eqnarray*}
for $|z|\leqslant 4 x_n$.
It follows that
\begin{eqnarray}\label{250}
\left| \sigma_n^2(z)-1\right| \leqslant C_{10}
\frac{g_n^2(8x_n)+g_n(8x_n)}{ x_{n}^2 (1-g_n(8x_n))}
\end{eqnarray}
for $|z|\leqslant 4 x_n$.

Let $\{h_n\}$ be a sequence of positive numbers.
Let $\overline{S}_n$ be a random variable conjugate to $S_n/\sqrt{B_n}$, i.e.
$\overline{S}_n$ has the following distribution function
$$ \mathbf{P} \left( \overline{S}_n < x\right) = 
\frac{1}{\varphi_n(h_n)} \int\limits_{-\infty}^x  e^{h_n u} d \mathbf{P}(S_n<u \sqrt{B_n}),
\quad x \in \mathbb{R}.
$$

Note that $\mathbf{E} \overline{S}_n = m_n(h_n)$ and
$D \overline{S}_n = \sigma^2_n(h_n)$. In the sequel, we take 
$h_n$ such that relations (\ref{240}) и (\ref{250}) will yield
$m_n(h_n)= h_n + o(1)$ and $\sigma^2_n(h_n)=1+o(1)$. 
Hence, we investigate the distance between the standard normal law
and the distribution of $ \overline{S}_n $ centered at and
normalized by main terms of the mean and the variance.
Denote
$$ R_n(v) = \mathbf{P}\left(\overline{S}_n-h_n<v\right)-\Phi(v), \quad v \in \mathbb{R},
$$
and estimate
$$ \Delta_n = \sup\limits_{v \in \mathbb{R}} | R_n(v)|.$$

Put
$$ \psi_n(z) = \mathbf{E} e^{ z (\overline{S}_n-h_n)}
= e^{- z h_n}
\frac{\varphi_n\left( z +h_n\right)}{\varphi_n(h_n)}, \quad z \in \mathbb{C}.
$$
It is clear that $\psi_n(it)$ is a c.f. of the random variable $ \overline{S}_n-m_n(h_n) $.

We have 
\begin{eqnarray*}
&&
\left| e^{-z^2/2} \psi_n(z) - 1 \right| = 
\frac{e^{h_n^2/2}}{\varphi_n(h_n)}
\left|
e^{-(z+h_n)^2/2} \varphi_n(z+h_n) - e^{-h_n^2/2} \varphi_n(h_n)
\right|
\\ &&
= \frac{e^{h_n^2/2}}{\varphi_n(h_n)}
\left| f_n(z+h_n)-f_n(h_n)\right|
\leqslant
\frac{e^{h_n^2/2}}{\varphi_n(h_n)}
\sum\limits_{k=1}^\infty |a_{nk}| \left| (z+h_n)^k-h_n^k \right|
\end{eqnarray*}
for $|z|+h_n\leqslant 4 x_n$.   Since
\begin{eqnarray*}
&& (z+h_n)^k - h_n^k = 
z (z+h_n)^{k-1} + h_n (z+h_n)^{k-1}  - h_n^k \\ &&
=z (z+h_n)^{k-1} + z h_n (z+h_n)^{k-2}  + h_n^2 (z+h_n)^{k-2}  - h_n^k 
\\ &&
=\cdots = 
z \sum_{j=1}^k (z+h_n)^{k-j} h_n^{j-1},
\end{eqnarray*}
we obtain $ |(z+h_n)^k - h_n^k| \leqslant k |z| (4 x_n)^{k-1}$
for $|z|+h_n\leqslant 4 x_n$.
It follows from relations (\ref{210}) and (\ref{220}) that
\begin{eqnarray*}
\left| e^{-z^2/2} \psi_n(z) - 1 \right| \leqslant
\frac{e^{h_n^2/2}}{\varphi_n(h_n)} |z|
\sum\limits_{k=1}^\infty |a_{nk}|   k (4 x_n)^{k-1} 
\leqslant
\frac{1}{1-g_n(h_n)} |z| 
 C_8  \frac{g_n(8 x_n)}{ x_{n}}
\end{eqnarray*}
for $|z|+h_n\leqslant 4 x_n$. 
Putting $z=it$, we get
\begin{eqnarray*}
\left| \psi_n(it) - e^{-t^2/2} \right|
\leqslant
C_8 |t| e^{-t^2/2} \frac{g_n(8 x_n)}{ x_{n} (1-g_n(h_n))}
\end{eqnarray*}
for all $|t|\leqslant 2 x_n$ and $|h_n|\leqslant 2 x_n$.
By the Esseen inequality, we have
\begin{eqnarray}\label{200}
\Delta_n = \sup\limits_{v \in \mathbb{R}} | R_n(v)|
\leqslant \frac{1}{\pi} \int\limits_{-2 x_n}^{2 x_n } 
\left| \psi_n(it) - e^{-t^2/2} \right| \frac{dt}{t} + \frac{24}{(2 \pi)^{3/2} x_n}
\leqslant
C_{11} \frac{1}{x_n}.
\end{eqnarray}

Furthermore,
\begin{eqnarray}
&& \nonumber
\mathbf{P}\left(S_n \geqslant m_n(h_n)\sqrt{B_n}\right) =
\varphi_n(h_n) 
\int\limits_{m_n(h_n)}^\infty  e^{-h_n u} d \mathbf{P}(\overline{S}_n<u)
\\ && \nonumber
=
\varphi_n(h_n) e^{- h_n^2}
\int\limits_{m_n(h_n)-h_n}^\infty  e^{-h_n v } 
d \mathbf{P}\left(\overline{S}_n-h_n<v\right)
\\ && \label{215}
=
\varphi_n(h_n) e^{- h_n^2}
\int\limits_{m_n(h_n)-h_n}^\infty  e^{-h_n v } 
d \left(\Phi(v) + R_n(v)\right).
\end{eqnarray}

We have
\begin{eqnarray} 
&& \nonumber
\int\limits_{m_n(h_n)-h_n}^\infty \!\!\!\!\!\! e^{-h_n v } d \Phi(v) 
= \frac{e^{h_n^2/2}}{\sqrt{2\pi}} 
\int\limits_{m_n(h_n)-h_n} \!\!\!\!\!\!
e^{-(v+h_n)^2/2} dv
\\ && \label{225}
=\frac{e^{h_n^2/2}}{\sqrt{2\pi}} \int\limits_{m_n(h_n)}^\infty\!\!\!\!\!\!
e^{-v^2/2} dv
\sim \frac{e^{(h_n^2-m_n^2(h_n))/2}}{\sqrt{2\pi} m_n(h_n)} ,
\end{eqnarray}
provided $m_n(h_n)\to\infty$. Moreover,
\begin{eqnarray}\label{235}
\left|\int\limits_{m_n(h_n)-h_n}^\infty\!\!\!\!\!\!  e^{-h_n v } d R_n(v)\right| =
\left| R_n(m_n(h_n)-h_n) -\!\!\!\! 
\int\limits_{m_n(h_n)-h_n}^\infty \!\!\!\!\!\! R_n(v) d (e^{-h_n v })\right|
\leqslant 2 \Delta_n.
\end{eqnarray}

Put $x_n=u_n \varrho_n$, where $\varrho_n\to\infty$ enough slowly to satisfy
$x_n\to \infty$, $x_n^3=o(\sqrt{n}/\gamma_n)$ and $x_n=o(\sqrt{B_n}/M_n)$.
Note that in view of relation (\ref{210}), we have $g_n(8 x_n) = o(1)$.
Let $h_n$ be a solution of the equation
\begin{eqnarray}\label{eqn}
m_n(h_n) = u_n.
\end{eqnarray}
The function $m_n(h)$ is strictly increasing, $m_n(0)=0$ and, by relations (\ref{240}) and (\ref{210}), 
the inequalities
$m_n(4 x_n) = 4 x_n + o(x_n^{-1}) \geqslant 2 x_n > u_n$ hold for
all sufficiently large $n$. It follows that the unique solution of equation (\ref{eqn}) exists
for all sufficiently large $n$. Moreover, relation (\ref{240}) yields that
$$ h_n = u_n + o(x_n^{-1})
$$
and
$$ m_n^2(h_n) -h_n^2 = (h_n+o(x_n^{-1}))^2 - h_n^2 = o(h_n x_n^{-1}) =o(1).
$$
It follows from relations (\ref{200})---(\ref{235}) and (\ref{210}) that
\begin{eqnarray*}
&& \nonumber
\mathbf{P}\left(S_n \geqslant u_n\sqrt{B_n}\right) =
\varphi_n(h_n) e^{- h_n^2} \left(  \frac{e^{(h_n^2-m_n^2(h_n))/2}}{\sqrt{2\pi} m_n(h_n)} (1+o(1))
+ O(\Delta_n)\right)
\\ &&
= e^{-\frac{h_n^2}{2}}
\left(\frac{1+o(1)}{\sqrt{2\pi} m_n(h_n)}   +O(x_n^{-1}) 
\right)(1+o(1))
\\ &&
= e^{- \frac{u_n^2}{2} +o(1)}
\left(\frac{1+o(1)}{\sqrt{2\pi} u_n}   +o(u_n^{-1}) 
\right)(1+o(1)) = (1-\Phi(u_n))(1+o(1)).
\end{eqnarray*}

Theorem \ref{t1} is proved. $\Box$

Finally, we mention some unsolved problems. In Frolov, Martikainen and Steine\-bach [18],
one can find more exact results on large deviations for sums of independent random variables
in the scheme of series. In there, the conditions are imposed on the logarithms of m.g.f.'s 
of summands. Now we can not adopt the techniques from there to combinatorial sums. We see from 
relation (\ref{efi}) that the m.g.f. of $S_n/\sqrt{B_n}$ is the permanent of the
matrix $\|\mathbf{E} e^{z X_{nij}/\sqrt{B_n}}\|$. Above, the method of the investigation 
of the behaviour for this permanent implied bounds with $ \gamma_n/\sqrt{n} $ instead of analogues
of the Lyaponov ratios. The second problem is that the proof in [18] involves some bounds
in CLT which variants for combinatorial sums are unknown. Solutions of these
problems could yield more exact results under weaker conditions.

\bigskip

{\bf References}

{\parindent0mm
{\begin{itemize}
\footnotesize
\item[{[1]}]
Wald A., Wolfowitz J.,1944. Statistical tests based on permutations of
observations. Ann. Math. Statist. 15, 358--372.
\item[{[2]}]
Noether G.E., 1949. On a theorem by Wald and  Wolfowitz.
 Ann. Math. Statist. 20, 455--458.
\item[{[3]}]
Hoeffding W., 1951. A combinatorial central limit theorem. Ann. Math. Statist. 22,
558--566.
\item[{[4]}]
Motoo M., 1957. On Hoeffding's combinatorial central limit theorem.
Ann. Inst. Statist. Math. 8, 145--154.
\item[{[5]}]
Kolchin V.F., Chistyakov V.P. (1973) On a combinatorial limit theorem.
Theor. Probab. Appl. 18, 728-739.
\item[{[6]}]
Bolthausen E., 1984. An estimate of the remainder in
a combinatorial central limit theorem.
Z. Wahrsch. verw. Geb. 66, 379--386.
\item[{[7]}]
von Bahr B., 1976. Remainder term estimate in a
combinatorial central limit theorem.
Z. Wahrsch. verw. Geb. 35, 131-139.
\item[{[8]}]
Ho S.T., Chen L.H.Y., 1978. An $L_p$ bounds for the remainder in
a combinatorial central limit theorem.
Ann. Probab. 6, 231--249.
\item[{[9]}]
Goldstein L., 2005. Berry-Esseen bounds for combinatorial central limit theorems
and pattern occurrences, using zero and size biasing. J. Appl. Probab. 42, 661--683.
\item[{[10]}]
Neammanee K., Suntornchost J., 2005. A uniform bound on
a combinatorial central limit theorem.
Stoch. Anal. Appl. 3, 559-578.
\item[{[11]}]
Neammanee K., Rattanawong P., 2009. A constant on a uniform bound of
a combinatorial central limit theorem. J. Math. Research 1, 91-103.
\item[{[12]}]
Chen L.H.Y., Goldstein L., Shao Q.M., 2011. Normal approximation
by Stein's method. Springer.
\item[{[13]}]
Chen L.H.Y., Fang X. (2015) 0n the error bound in a
combinatorial central limit theorem. 
Bernoulli, 21, N.1, 335-359.
\item[{[14]}]
Frolov A.N., 2014. Esseen type bounds of the remainder in a combinatorial CLT.
J. Statist. Planning and Inference, 149, 90--97.
\item[{[15]}]
Frolov A.N. (2015a) Bounds of the remainder in a combinatorial central limit theorem.
Statist. Probab. Letters 105, 37-46.
\item[{[16]}]
Frolov A.N. (2015b) On the probabilities of moderate deviations for combinatorial sums.
Vestnik St. Petersburg University. Mathematics, 48, No. 1,  23-28. Allerton Press, Inc., 2015.
\item[{[17]}]
Frolov A.N. (2017) On Esseen type inequalities for combinatorial random sums. Communications in Statistics -Theory and Methods. 46 (12), 5932-5940. 
\item[{[18]}]
Frolov A.N., Martikainen A.I., Steinebach J. (1997) Erd\"{o}s--R\'{e}nyi--Shepp type laws in non-i.i.d.
case. Studia Sci. Math. Hungar.  34, 165--181.
\end{itemize}
}
}
\end{document}